# Is mathematics invading human cells?
*Impressions from a collaboration with diabetes doctors*
**Bernhelm Booß-Bavnbek, Roskilde University, Denmark**

Jointly with another mathematician, a biophysicist and two diabetes doctors, the author has just released a textbook, *BetaSys - Systems Biology of Regulated Exocytosis in Pancreatic ß-Cells*[1], in which a broad international team summarizes the state of our current understanding of the cell-physiological events accompanying both successful – and impaired insulin secretion. *The Mathematical Intelligencer* has asked Bernhelm Booß-Bavnbek to describe some of his experiences as a mathematician cooperating with diabetes specialists.

## Advanced equipment and basic ignorance

Along with space exploration and military and civilian nuclear power design, medical devices belong to the mathematically most sophisticated areas of modern technology. There are, for example, many mathematicians that have either contributed or could have contributed to magnetic spin resonance imaging (MRI), and there is hardly a single mathematician who masters all the math involved in that technology. The same goes for electron tomography, multi-beam confocal laser microscopy and many other advanced devices. Medicine has become a mathematical discipline. The ominous military-industrial complex has grown into an eminent mathematical sick-and-health industry.

But mathematics is encapsulated in the apparatus. Whether it is about a specific diagnosis or treatment, most patients will probably, at least when you come from mathematical physics, be surprised at how little medical science seems to really know and understand about the individual diseases. It is quite normal that a doctor must simply experiment - or just stick to an established symptom diagnosis and symptom treatment. Without a detailed identification of the real causes of the individual patient's ailment, often a successful treatment, defined as cure, is unattainable.

Physics can also be complicated and in many cases without established answers. But in physics there is after all only a very short list of "First Principles", one must stick to. There we have relatively well-defined interfaces between certain knowledge, reasoned or vague presumption and ignorance. And in most cases, our ignorance in physics can be condensed in some mathematical equations (which we, however, may not immediately fully understand). This is not the case so in medicine.

## Specific nature of the use of mathematics in cell research

*The strong medical pull*. From pure mathematical research, we know the feeling of being pulled forward by an overarching issue: the relationship between local and global properties, between the smooth and the continuous, between analytic and algebraic methods, the Four Color Problem, the Poincaré Conjecture, the Riemann Hypothesis, the Clay Millennium Problems. Of course, we would never admit such personal ambitions in public. But for me there is no doubt about the role that major well-stated problems play and have played in the design of the career paths of many mathematicians, at least indirectly and in daydreams: with many doubts and a persistent feeling of self-deception and of fighting against mountains - or windmills.

---

[1] Booß-Bavnbek, B.; Klösgen, B.; Larsen, J.; Pociot, F.; Renström, E. (eds.), *BetaSys - Systems Biology of Regulated Exocytosis in ß-Cells*, series: Systems Biology, Springer, Berlin-Heidelberg-New York, 2011, XVIII, 558 pages, 104 illustr., 53 in color. With online videos and updates. ISBN 978-1-4419-6955-2. [1] Comprehensive review in *Diabetologia*, DOI 10.1007/s00125-011-2269-3.

Working as a mathematician with diabetes doctors is different. A bristling cascade of medical issues pulls the research forward: For nearly 90 years we have known that lack of secretion of the hormone insulin is one of the many serious issues in both diabetes type 1 (juvenile) and type 2 (obesity and age driven). For a large group of these patients, actually insulin is produced in pancreatic β-cells and stored in thousands of mini bags, *vesicles*, in the cell's interior. But the cells do not respond correctly to external stimuli with the actual secretion, called *regulated exocytosis*. That manifests itself in elevated blood sugar, which can be tasted and measured by urine sample. That has now been done for nearly four thousand years.[2] We call it a *symptom diagnosis* because the diagnosis says nothing about the wide range of causes which may underlie the lack of uptake of glucose in the muscles.

Previously, failure of insulin secretion automatically led to weakening the muscles, inflammation of the extremities, loss of vision and the body's final decay. Since the discovery of insulin, this tragic development can be countered by artificial supply of insulin by injection several times a day. We call it a *symptom treatment* because it is not even attempted to cure the patient or to make an effort to restore the body's own insulin secretion. Some claim that the relative success of the overall symptom diagnosis and symptomatic treatment of diabetes has blocked patient-centered, individualized diagnosis and treatment.

In any case collaboration with diabetes doctors is a powerful experience for a mathematician of continually being pulled forward by well-defined medical problems. Here it is simply to detect the functioning and system behavior of the regulated exocytosis in healthy β-cells and to identify everything that can stand in the way in the case of weakened β-cells, see Fig. 1. The purpose is clear: mathematician, please come and help find the way to an earlier and more specific diagnosis and a cure or alleviation of the specific failure!

*The technological push*. The technological push is not completely unfamiliar in mathematics, we may think: readily available electronic journals, large user-friendly collections of mathematical preprints and reviews, efficient numerical software packages, homemade LaTeX editing can put us under pressure as mathematicians. But it's nothing compared to the immense technological pressure cell research is subject to: with each new generation of equipment, oceans of new data inundate on quite different length scales. Rapidly expanding technology-pushed innovations are, e.g., about individual genes in the DNA, about proteins and about electrical cell membrane processes, but also the structure and function of a β-cell as a whole are attempted to be described in momentary images (by electron tomography) or dynamic sequences (by tracking of properly primed nanoparticles in living cells).

*Heavy preponderance of ad-hoc perceptions*. There is no shortage of heroic attempts by some scientists to bring a little order and overview into this real wild jungle of data. Most tries, however, restrict themselves to ad-hoc fiddled perceptions of unconfined creativity à la: "it should probably be the cell nucleus that controls the process" or "there is a certain rate, which determines the transition between one stage and the next" or "a correlation between the one process and another process is unquestionable". Explanations hold until overtaken by new data and will then be "adjusted". They will never be falsified because they are freestanding and variable and not, as we are accustomed to from the world of physics, tied by head and limbs to the basic physical laws and the geometric properties of the three-dimensional space. The only quality criterion is whether a model *looks like* the known observations or can be tuned to coincide with them. It is a free kingdom of modeling, admitting fancied ghosts to explain actual observations, but a nightmare when looking for descriptions and explanations, offering a certain shelf life and a theoretical check for errors.

---

[2] The earliest preserved report (in Bendex Ebbell's Copenhagen interpretation of 1937) is from the Egypt *Ebers Papyrus* of 1536 BCE., instruction 197, column 39, line 7, reproduced in all its ambiguity on http://biology.bard.edu/ferguson/course/bio407/Carpenter_et_al_(1998).pdf.

*The phylogenetic heritage*. Our insulin-producing β-cells are among the most differentiated human cells. They are closely packed with a zoo of different types of organelles. Insulin-like peptides can be detected in our distant invertebrate ancestors who have been around for more than 600 million years. Something resembling pancreases with a kind of insulin-producing β-cells already exist in the hagfish, which have existed for more than 500 million years. For every discovery it must be feared that a new observed process, a new measured quantity is quite irrelevant. Maybe you just hit a relic, a ruin of the historical development, which has no importance anymore. Of course, this type of confusion was also met in the history of physics. How long has it taken to assign to meteors and comets their place in our conception of the solar system or to remove Pluto from the list of planets formed in our solar system? However, while the ruins and relics ideally sharpen the mind in simple research fields such as physics and astronomy, they can be extremely confusing and even completely block medical research. Again and again one senses that we mathematicians coming from the outside possibly are too early. Perhaps we had better wait for another 150 or 200 years until the research has separated essential processes from nonessential processes, before we at last can begin the serious work.

*Lack of universality*. What strikes me most in mathematical cell physiology is the lack of any *universality* or *scale invariance*. In the world of physics, Maxwell's equations apply both for high-frequency radio waves and low frequency voltage in power plants; the Navier-Stokes equations apply for both the continental atmospheric phenomena and the whirling around a ship hull. In physics, we have a field concept relating point measurements with spatially widespread events. It is not (yet?) so in cell physiology.

Of course there are cross connections between what we know about β-cell function and our genetic data, our conceptions of the mode of operation of single organs (like the pancreas) and a body's, an organism's behavior and the performance of a whole population. For example, genetic data are just collected by epidemiological studies of large populations and the feedback is well studied between nutrient intake, liver and brain response and the secretion signaling. But - apart from the universality of the applied statistical methods for parameter estimation and hypothesis testing – all the met methods are closely tied to a specific biological level, a particular length and time scale. We know such a hopeless situation also from mathematical physics with the seeming incompatibility between the mathematical theories of gravitation and quantum mechanics. That might be considered as wounds in physics, but it is a unique wound. In diabetes research, we have hundreds of such cracks and ditches where no one knows if there is a bridge or how it then would be built.

*Volatility*. Medical biology, as it is conducted today is a huge undertaking with a myriad of articles published every year. Not many of them will be quoted after two years. That's probably the reason that a key parameter for bibliometric research information, the *impact factor*, only examines the current references to papers that are not more than just these two years old. Sure enough, the overall goal, the understanding of life and death, of health and illness, is long-lasting. But the angles of attack change constantly and appear frequently as dictated by some observational techniques that have just now come to use. The subject seems to be characterized by the absence of established and general traditions. As practiced today, cell physiology is a young subject which is just establishing itself. Accidental discoveries seem to play a major role. We recognize that also from physics, where, e.g., the discovery of high-temperature super-conductivity in conventional insulating ceramic materials by Bednorz and Müller in 1986 could hardly be characterized as the result of deep theoretical considerations. However, random breakthroughs occur without doubt more often in biomedicine.

*Systems thinking versus reductionism*. It goes without saying that a strictly reductionist program is needed in medical research, if the current packing of medical ignorance in ad-hoc assumptions shall be replaced by falsifiable references to basic physical laws. But I must also acknowledge that most bodily functions and processes involve many different cell components, neighboring cells, various organs and the whole organism in an interaction. Understandably, the holistic slogan of *systems biology* has become popular, and great expectations are attached to it.

Both programs will reveal exciting new facts and relations. Both approaches offer the mathematician rich working opportunities. To me, the most promising direction is somewhere in the middle: maybe a *focused systems biology* will show its ability to touch the wall, knock a hole in it, and achieve a breakthrough. That has not happened yet. The hope is to develop a medicine and a biology that simplifies in a reductionist way; fearlessly ignores some probably relevant aspects; and focuses on a limited range of processes; but in turn lets itself holistically and equally fearlessly be confronted with a multitude of levels and a diversity of length and time scales all at once.

## Mathematical helping hand

What place ‚then, has a mathematician in this environment?

*The daily practice*. Just as in engineering, economics or anywhere else, also in cell physiology the daily mathematical exercise consists of the estimation of some parameters, testing the significance of some hypotheses and designing compartment models for the dynamics of coupled quantitative variables. Often, the role of mathematics is to check whether a random discovery delivers what it promised.

Numerical problems can rapidly pile up when one wants, e.g., to simulate a fusion process of a simple insulin vesicle to the plasma membrane of the β-cell throughout the process: the bending of the plasma membrane into a dimple, the coupling of the vesicle to the dimple, the coalescence of vesicle and plasma membrane during the hemifusion, the formation of the fusion pore for emitting the insulin molecules, and dissolving the vesicle remains into the plasma membrane. The reason for the numerical problems is that we are at a mesoscale: the characteristic lengths vary from 1 nm for the lipid heads, to 7nm for the strength of the membrane bilayers of lipids - to 100-250 nm for the insulin vesicle diameters. Thus, the relevant lengths of regulated exocytosis considerably exceed the lengths that chemists have mastered using *Molecular Dynamics* (MD). It is even worse with the time scale, because a simple β-cell responds to glucose stimulation by insulin secretion over 25-30 minutes. And everything is in three dimensions, see Fig. 2. This requires the development of special software to aggregate both space and time intervals to something that existing current computers can work with.[3]

*The dual role of mathematicians*. A mathematician coming from the outside must be humble in front of the immense calibration and programming work that underlies such models. It's hard not to succumb to the fascination of the "lively" graphical output of such simulations. Respectfully and humbly, we should make our *tool box* available and fearlessly lend a hand when needed. But we must not abandon our *mathematical way of viewing*, our acquired competence to inquire into the basis for the modeling and the simulations. We must remain skeptical and question everything by cross-checking calculations; insisting on relating corresponding phenomena with each other; and using our imagination to devise quite simple physical mechanisms that have the ability to generate the world of complex phenomena that we observe.

---

[3] J. Shillcock, *Probing cellular dynamics with mesoscopic simulations*, in [1], pp. 459-473.

*The falsifying and heuristic function of mathematics*. There are many jokes about the sharp-nosed mathematicians who check up something and afterwards, sometimes annoyed, sometimes smiling note that the biologists' data and assumptions do not fit together. This gives mathematicians a reputation of pettiness and pedantry, but it is perhaps our most important contribution to all biomedical fields. With such a know-all tone, between 1616 and 1628, William Harvey falsified the prevalent notions about the cardio-vascular system and discovered the *arithmetic existence* of the blood capillaries that connect arteries and veins - 40 years before Marcello Malpighi's light microscope confirmed their *histological reality*.[4]

Similarly, e.g., a harmonic analysis[5] of observed electrical vibrations (calcium oscillations) in β-cells just before secretion indicates that these fluctuations are not only an expression of pulsing influx of calcium ions through the plasma membrane, but - contrary to popular perception - may also result from a pulsing violent "splashing" of these ions between the cell's internal calcium organelles such as mitochondria and the endoplasmic reticulum. Hence, a purely mathematical realization of an inconsistence can move the focus from, I must admit, more easily and directly measurable local electrical membrane processes (measurement of the change of the static potential over time using the patch clamp) to cell-internal global and long-range electro-dynamic processes (measurement of fluctuating magnetic field strengths) and give the exocytose research a new approach.[6]

*Model-based and simulated measurements*. Many biomedical quantities cannot be measured directly. That is due to the subject matter, here the nature of life, partly because most direct measurements will require some type of fixation, freezing and killing of the cells, partly due to the small length scale and the strong interaction between different components of the cell. Just as in physics since Galileo Galilei's determination of the fall law by calculating "backwards" from the inclined plan, one must also in cell physiology master the art of model-based experiment design. Let us, e.g., look at the eight to twelve thousand densely packed insulin vesicles in a single β-cell. They all must reach the plasma membrane within a maximum of 30 minutes after stimulation, to pour out their contents. Let us ignore the many processes taking place simultaneously in the cell and consider only the basic physical parameter for transport in liquids, namely the viscosity of the cell cytosol. From measurements of the tissue (consisting of dead cells) we know the magnitude of viscosity of the protoplasma, namely about 1 milli-pascal-seconds (mPa s), i.e., it is of the same magnitude as water at room temperature. But now we want to measure the viscosity in living cells: before and after stimulation; deep in the cell's interior and near the plasma membrane; for healthy

---

[4] For details, cf. the box *Harvey's arithmetical microscope* in J. T. Ottesen, *The mathematical microscope – making the inaccessible accessible*, in [1], pp. 97-118, here p. 99.
[5] L. E. Fridlyand and L. H. Philipson, *What drives calcium oscillations in β-cells? New tasks for cyclic analysis*, in [1], pp. 475-488.
[6] D. Apushkinskaya et al., *Geometric and electromagnetic aspects of fusion pore making*, in [1], pp. 505-538.

and stressed cells.

It serves no purpose to kill the cells and then extract their cytosol. We must carry out the investigation *in vivo* and *in loco*, by living cells and preferably in the organ where they are located. The medical question is clear. So is the appropriate technological approach, since techniques have been developed which allow iron oxide nanoparticles of a diameter up to 100 nm to be brought inside these most vulnerable β-cells without destroying them. It happens with a low frequency (around 10 Hz) electromagnetic dynamic field generator that makes nanoparticles so to speak, to "roll" on the surface of the cells until they hit a willing receptor and get approach to the cell interior across the plasma membrane. These particles are primed with appropriate antigens and with a selected color protein, so that their movements within the cell can be observed with a confocal multi-beam laser microscope which can produce up to 40 frames per second. The periods of observations are only relatively short, perhaps a maximum of 8-10 minutes - before these particles are captured by cell endosomes and delivered to the cells' lysosomes for destruction and consumption of their color proteins.[7]

The simplest mathematical method to determine the viscosity of the cytosol in vivo would be just to pull the magnetized particles with their fairly well-defined radius *a* with constant velocity *v* through the liquid and measure the applied electromagnetic force *F*. Then the viscosity η is obtained from *Stokes' Law* $F = 6\pi a \eta v$. The force and the speed must be small so as not to pull the particles out of the cell before the speed is measured and kept constant. Collisions with insulin vesicles and other organelles must be avoided. It can only be realized with a low-frequency alternating field. But then Stokes' Law must be rewritten for variable speed, and the mathematics begins to be advanced. In addition, at low-velocity we must correct for the spontaneous Brownian motion of particles. Everything can be done mathematically: writing the associated stochastic *Langevin equations* down and solve them analytically or approximate the solutions by Monte Carlo simulation.[8] However, we rapidly approach the equipment limitations, both regarding the laser microscope's resolution and lowest achievable frequency of the field generator.

So we might as well turn off the field generator and be content with intermittently recording the pure *Brownian motion* of a single nanoparticle in the cytosol! As shown in the two famous 1905/06 papers by Einstein[9], the motion's variance (the mean square displacement over a time interval of length τ) $\sigma^2 = <\mathbf{x}^2> = E(|\mathbf{x}(t_0+\tau) - \mathbf{x}(t_0)|^2)$ of a particle dissolved in a liquid of viscosity η is given by $\sigma^2 = 2D\tau$, where

$$D = \frac{k_B T}{6\pi a \eta}$$

denotes the diffusion coefficient with Boltzmann constant $k_B$, absolute temperature *T* and particle radius *a*. In statistical mechanics, one expects $10^{20}$ collisions per second between a single colloid of 1 μm diameter and the molecules of a liquid. For nanoparticles with a diameter of perhaps only 30 nm, we may expect only about $10^{17}$ collisions per second, still a figure so large as to preclude registration. There is simply no physical observable quantity $<\mathbf{x}^2>$ at the time scale $\tau = 10^{-17}$ seconds. But since the Brownian motion is a Wiener process with self-similarity we get approximately the same diffusion coefficient and viscosity estimate, if we, e.g., simply register 40 positions per second. Few measurements per second are enough. *Enough is enough*, we can explain to the experimentalist, if he/she constantly demands better and more expensive apparatus.

---

[7] Details will be disclosed in a US patent in preparation by M. Koch et al.
[8] F. Schwabl, *Statistical Mechanics*, Springer, Berlin-Heidelberg-New York, 2006; A. R. Leach, *Molecular Modelling - Principles and Applications*, Pearson Education Ltd., Harlow, 2001, Chapter 7.8.
[9] A. Einstein, 'Über die von der molekularkinetischen Theorie der Wärme geforderte Bewegung von in ruhenden Flüssigkeiten suspendierten Teilchen', *Ann. Phys.* **17** (1905) 549-561; 'Zur Theorie der Brownschen Bewegung', *Ann. Phys.* **19** (1906) 371-381. Both papers have been reprinted and translated several hundred times.

This is beautifully illustrated by a small MatLab program (see Box 1), which first generates a Wiener process with a given variance $\sigma^2$ and then estimates the variance from the zigzag curves generated by taking all points or every second or fourth. Note that $\sigma^2$ also can be estimated by the corresponding two-dimensional Wiener process of variance 3/2 $\sigma^2$, consisting of the 2D-projections of the three-dimensional orbits, as the experimental equipment also will do.[10]

Beautiful, but it is still insufficient for laboratory use: there we also must take into account the non-Newtonian character of the cytosol of β-cells. These cells are, as mentioned, densely packed with insulin vesicles and various organelles and structures. Since the electric charge of iron oxide particles is neutral, we can as a first approximation assume a purely elastic impact between particles and obstacles. It does not change the variance in special cases as M. Smoluchowski already figured out 100 years ago for strong rejection of particles by reflection at an infinite plane wall.[11] Here also computer simulations have their place to explore the impact of different repulsion and attraction mechanisms on the variance.

Now you can hardly bring just a single nanoparticle into a cell. There will always be many simultaneously. Thus it may be difficult or impossible to follow a single particle's zigzag path in a cloud of particles by intermittent observation. Also here, rigorous mathematical considerations may help, namely the estimation of the viscosity by a periodic counting of all particles in a specified "window".[12]

The goal of model-based measurements and computer simulations is both to obtain the desired quantity from available or realizable observations *and* to become familiar with the expected laboratory conditions. Calculations and simulations can make us on intimate terms with the expected results; can support the exploration of a range of a priori unknown conditions; and help to identify the best choice of free parameters such as particle diameter, temperature, area of focus etc.

*Need for new mathematical ideas*? I have described how important a wide solid mathematical competence is for success in everyday practice; both for the verification and falsification of current assumptions; and for model-based measurements and simulation. Overview and literature study are required, not originality, this *mother of banalities*, as a Ukrainian *bon mot* says.

But there is also a need for radically new mathematical ideas, especially ideas that can integrate the otherwise isolated and local observations and perceptions that characterize molecular biology. How localized events propagate from a position at the plasma membrane into a global process involving a myriad of ions, proteins and organelles far away and across the cell and let the essential event take place: the secretion, back at the plasma membrane? How does the communication taking place; the spread of a singularity; the amplification of a signal; and finally the creation of new forms? Many mathematical disciplines have something to offer, from algebraic geometry, stochastic processes and complex dynamics to parabolic and hyperbolic differential equations and free boundary value problems.[13]

---

[10] M. von Smoluchowski, 'Zur kinetischen Theorie der Brownschen Molekularbewegung und der Suspensionen', *Ann. Phys.* **21** (1906), 756-780, §9 gives – erroneously - the correction factor 4/π, i.e., the reciprocal value of the average shortening of a 3D-length in 2D-projection.
[11] M. von Smoluchowski, 'Einige Beispiele Brownscher Molekularbewegung unter Einfluß äußerer Kräfte', *Bull. Int. Acad. Sc. Cracovie, Mat.-naturw. Klasse A* (1913), 418-434.
[12] M. von Smoluchowski, 'Studien über Molekularstatisktik von Emulsionen und deren Zusammenhang mit der Brownschen Bewegung', *Sitzber. Kais. Akad. Wiss. Wien, Mat.-naturw. Klasse* **123**/IIa (Dec. 1914), 2381-2405. All three here cited papers by Smoluchowski are available on http://matwbn.icm.edu.pl/spis.php?wyd=4&jez=en.

## Conclusion

*How deep is the gap between mathematics and medicine?* Most mathematicians who have tried to work with doctors will confirm that cooperation is fairly smooth. You soon find a common language and common understanding in spite of widely different backgrounds. Understandably, one should not and cannot overstretch the patience of a clinical physician who has his or her patient here and now.

The relationship between mathematics and medicine has been quite tumultuous in the history of science. Important mathematicians and physicists as R. Descartes, D. Bernoulli, J. d'Alembert, H. Helmholtz, E. Schrödinger, I. Gelfand, and R. Thom have been attracted to biomedical questions and observations, but have also expressed their reservations. Important doctors, for examples one need only go through the list of Nobel Laureates, have apparently not suffered from math phobia, but rather retained a lifelong fondness for mathematical ideas and ways of seeing.

Maybe this understanding between physicians and mathematicians has deep roots in the past: counting and healing was, by all accounts, magicians' and medicine men's mysterious privilege in pre-scientific cultures. Both subjects were, however, unlike the previous conjuration spirit and belief in magic and the good or evil ghosts, carried by the same rationalistic spirit throughout Greek and Roman antiquity (perhaps beside the Asclepiades). Geometric and arithmetic ratios should be explained and not adored or cursed! In the same mind, Greek medicine has established itself as a strictly materialistic subject who described the disease course in purely objective, observable terms, and also envisioned solely objective reasons and pure physical treatment.[14]

*Tasks for mathematics education.* All higher educational institutions within mathematics have over the last years experiences that often more than half of their graduates were employed in the financial sector, especially to the mathematically delicate evaluation of options and other derivatives. Some university teachers have been just as pleased as their students over these quick appointments. Some went so far as to point to this new job market as an argument to attract new math students to their universities.

I agree with the series of critical contributions in the *Mathematical Intelligencer*: there is no reason to be proud to have trained some of our best students just to that task.[15] One alternative is to train

---

[13] For the last mentioned approach cf. D. Apushkinskaya et al. 2012, loc. cit. For a more fundamental approach to the geometry of biological amplification processes see also M. Gromov's many related and quite varied contributions of the last decade.

[14] Paul Diepgen, *Geschichte der Medizin. Die historische Entwicklung der Heilkunde und des ärztlichen Lebens*, Vol. 1, Walter de Gruyter & Co., Berlin, 1949, pp. 67-158; Fridolf Kudlien, *Der Beginn des medizinischen Denkens bei den Griechen, Von Homer bis Hippokrates*, Artemis, Zürich and Stuttgart, 1967; Fritz Jürss, *Geschichte des wissenschaftlichen Denkens im Altertum*, Akademie-Verlag, Berlin, 1982.

[15] M. Rogalski, 'Mathematics and finance: An ethical malaise', *Mathematical Intelligencer* **32**/2 (2010), 6-8; I. Ekeland, 'Response to Rogalski', *Mathematical Intelligencer* **32**/2 (2010), 9-10; J. Korman, 'Finance and mathematics: A lack of

our students in pure mathematics at its best. Perhaps an even better alternative is to direct students' attention to the many fascinating possibilities of cooperation in the medical world at the population level, e.g., in the study of infectious diseases and antibiotic resistant bacteria; at the organism and organ level, e.g., in the study of cardiovascular diseases; or at the cellular level, e.g., in the study of β-cells and other highly differentiated cell types.

## Thanks

I thank the doctors Hans-Georg Mannherz (Universität Bochum und Max-Planck-Institute Dortmund), Pierre de Meyts (NovoNordisk-Hagedorn Laboratories, Gentofte), Flemming Pociot (Glostrup Hospital, Copenhagen) and Erik Renström (University Hospital, Malmö) for many years of inspiration - and patience with me as a novice. Engineer Martin Koch has awakened my interest in the area and introduced me to the literature and to the electro-dynamic laboratory techniques. My Roskilde colleague Nick Bailey has done his to reduce the number of linguistic errors and ambiguities in this report.

## Illustrations:

**Fig. 1**. Two-phase secretion of insulin with three different β-cell modes. The figure shows at the bottom three β-cells in three different states. The smaller circles symbolize insulin vesicles. The graph on top shows the insulin secretion over time for a single cell. As the graph shows, the insulin secretion is explosive in the short first phase (mode i). In the longer second phase (mode iii), the secretion is rather constant and more evenly distributed. Between the two phases is the waiting state ii. As depicted in the β-cells at the bottom of the figure, the three molecular states are similar to each other. Consequently, they do not explain the order in the sequence of phases. It is that order which one now seeks to explain by means of mathematical models that involve the interplay between all processes. After Renström (2011) in [1], p. 40

**Fig. 2**. Cross section of the fusion of a vesicle of 28 nm diameter and a 100 x 100 $nm^2$ plane lipid bilayer in computer simulation, instantaneous snapshot 300 ns after the first vesicle-membrane contact. The simulation of a complete fusion event requires 4 CPU days on a processor. From Shillcock (2011) in [1], p. 468

**Box 1**. Stump of a MatLab program to visualize the self-similarity of a Wiener process, to be built upon for more realistic and investigative purposes. Modified from Gyurov, A., and R. Tokin, *Modeling the measurements of cytosol viscosity of pancreatic beta-cells by nanoparticle tracing in vivo*, Basic Studies in Natural Sciences, RU, Spring 2011.

---

debate', *Mathematical Intelligencer* **33**/2 (2011), 4-6. Related questions have been addressed in *SIAM-News* and *Mitt. Deutsch. Math.-verein*.

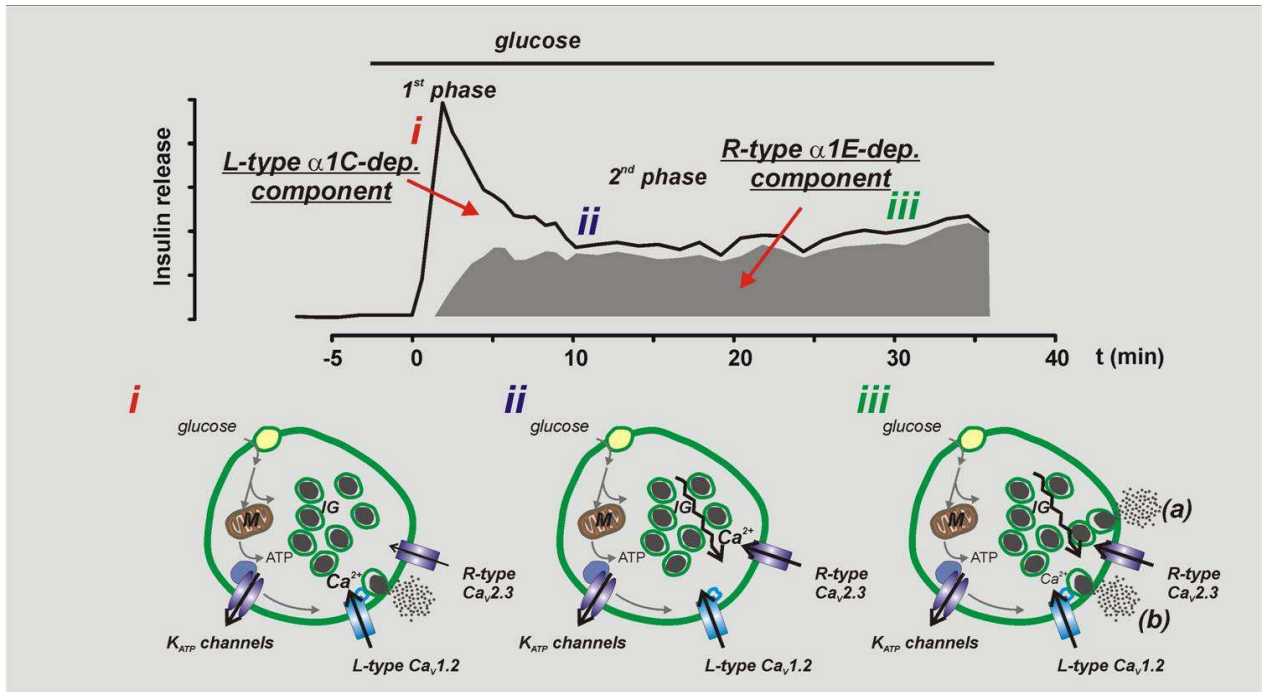

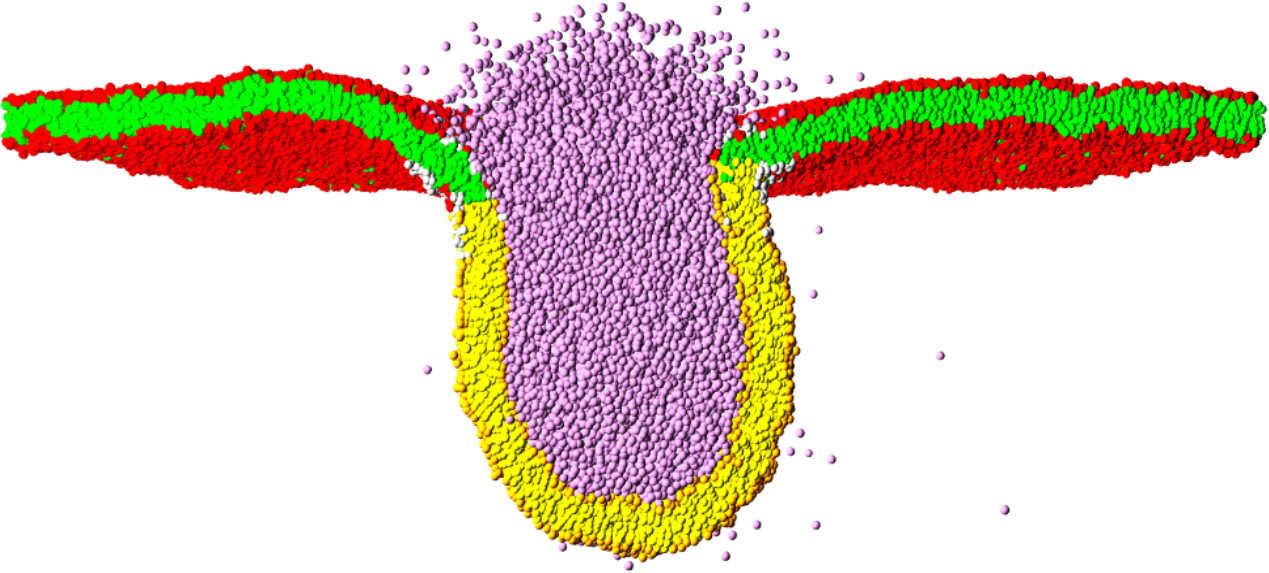

```matlab
% code1_2Dto3D_v2, , ag&rt, RUC NatBas, modified 8 August, 2011, bbb & jørgen larsen
% This program explores the most primitive way of estimating the viscosity of a fluid,
% e.g., the cytosol of pancreatic beta-cells in vivo, by tracing the Brownian motion {P(i)}
% of a single suspended nanoparticle by intermitting laser microscope registration.
%
% Problem: by computer simulation, to determine the expected precision of viscosity estimation in laboratory,
% depending (i) on the time resolution of the applied laser microscope, (ii) on the
% duration of possible observation, (iii) on the temperature: viscosity eta_temp = A exp(b/temp)
%
% Program design: Block I specifies laboratory parameters and calculates the corresponding
% diffusion constant D. Block II generates a sequence of Brownian 3D-motions, based on the
% calculated diffusion constant and a normal random generator. Each Brownian motion comes
% in three versions: all 40 frames per second, half of the frames, and a quarter of the
% frames. Block III calculates the observed square displacement between consecutive
% positions; derives estimates of D, and determines the viscosity eta on that basis. Block IV provides
% a graphical illustration of the self-similarity and of the considerable average maximal displacement
%
% Note: "clear" the command window before each run.
%% Block I: Input data---------------------------------------------
M=100;                          % number of tests of estimation quality
obs=4*60;                       % basic time length of a single video clip [sec]
N=obs*40;                       % number of measurements in high resolution
h=1/40;                         % time step in high resolution [sec]
t=(0:h:obs);                    % used for various resolutions
kb=1.38*10^(-23);               % Boltzmann constant
temp=310;                       % temperature [degrees K] = 37 [C]
visc=1;                         % viscosity [mPa s]
npr=0.02*10^(-6);               % NP radius [m]
D=1000*kb*temp/(6*pi*npr*visc)  % diffusion constant; recall the Einstein relation:
% 2*D*tau=MSD_tau=E(|P(i+tau)-P(i)|^2) for each time interval tau
% Input MSD_h=E(|P(1)-P(0)|^2) mean square displacement per time unit h
MSD_h=2*D*h                     % 2 * diffusion constant * unit time; OBS: h = 1/40
sigma_3D=sqrt((MSD_h)/3);       % factor for randn in coordinates to generate
% three-dimensional Wiener process of variance MSD_1
%%
% Block II: Generation of Brownian path in 3D with given variance sigma^2=MSD_h=2*D*h
for j=1:M
x(1)=0.0; y(1)=0.0; z(1)=0.0;   % Initial position of the Brownian nanoparticles
  for i=1:N                     % generating NP positions
  x(i+1)=x(i)+sigma_3D*randn; y(i+1)=y(i)+sigma_3D*randn; z(i+1)=z(i)+sigma_3D*randn;
  L(i)=norm([x(i+1)-x(i) z(i+1)-z(i)]); % array of displacements in x-z plane
  end
% 1/2th of the initial resolution
x2=x(1:2:end); y2=y(1:2:end); z2=z(1:2:end);
  for k=1:N/2        % generating vector with 1/2 jump-lengths
  L2(k)=norm([x2(k+1)-x2(k) z2(k+1)-z2(k)]);
  end
% 1/4th of the initial resolution
x4=x2(1:2:end); y4=y2(1:2:end); z4=z2(1:2:end);
  for m=1:N/4        % generating vector with 1/4 jump-lengths
  L4(m)=norm([x4(m+1)-x4(m) z4(m+1)-z4(m)]);
  end
% Block III: MSD observed in x-z plane, to be corrected by factor 3/2
MSDout=sum(L.^2)/length(L); MSD2out=sum(L2.^2)/length(L2);
MSD4out=sum(L4.^2)/length(L4);
% Estimated diffusion constant for different time resolution
D_obs=(3/2)*MSDout/(2*h);       % 40 slides/s
D2_obs=(3/2)*MSD2out/(2*2*h);   % 20 slides/s
D4_obs=(3/2)*MSD4out/(2*4*h);   % 10 slides/s
coef=1000*kb*temp/(6*pi*npr);
% Estimated viscosity for total observation time obs
% based on [40 20 10]  slides/s
visc_est(j)=coef/D_obs; visc2_est(j)=coef/D2_obs; visc4_est(j)=coef/D4_obs;
end
% Predicted error of viscosity estimates
visc_StandD=[sqrt(mean((visc_est-visc).^2)) sqrt(mean((visc2_est-visc).^2)) sqrt(mean((visc4_est-visc).^2))]
% Typical results for naive estimation, based only on the analysis of subsequent pairs of positions
% obs=1 sec: visc_variance = [0.1306 0.1611 0.2739] useless for all resolutions
% obs=10 sec: visc_variance = [0.0491 0.0738 0.0883] useless, expected error too large for all time resolutions
% obs=60 sec: visc_variance = [0.0184 0.0265 0.0380] still useless for wanted distinctions \pm few percent
% obs=4*60 sec: visc_variance = [0.0093 0.0137 0.0205] useless for low time resolution
% obs=10*60 sec: visc_variance = [0.0078 0.0081 0.0128] yields meaningful viscosity estimation
% even for low time resolution, i.e., greater obs can compensate for lower resolution.
% Block IV: Illustration of self-similarity for [40 20 10] slides/s, best for obs=1
hold on
plot3(x,y,z) plot3(x2,y2,z2,'r') plot3(x4,y4,z4,'g') xlabel('x') ylabel('y') zlabel('z')
```